\documentclass[12pt]{article}
\usepackage{epsfig}
\usepackage{amsmath,amssymb,amsthm,amscd}
\newcommand{\IP}{\relax{\rm I\kern-.18em P}}
\newcommand{\IR}{\relax{\rm I\kern-.18em R}}

\newcommand{\Cal}{{\cal }}
\newcommand{\Ž}{{\' e}}
\newcommand{\}{{\`e}}
\newcommand{\ˆ}{{\`a}}

\newcommand{\}{{\`u}}
 \newcommand{\Pic}{{\rm Pic}}
 \newcommand{\Isom}{{\rm Isom}}
 \newcommand{\Aut}{{\rm Aut}}
 \newcommand{\GL}{{\bf GL}}
 \newcommand{\PGL}{{\bf PGL}}
 \newcommand{\Spec}{{\rm Spec}}
 \newcommand{\pgcd}{{\rm pgcd}}
 \newcommand{\Ext}{{\rm Ext}}
 \newcommand{\Res}{{\rm Res}}
 \newcommand{\Div}{{\rm Div}}
\newcommand{\Bb}{{\bf B}}
\newcommand{\R}{{\rm R}}

\newcommand{\be}{\begin{equation}}
\newcommand{\ee}{\end{equation}}
\newcommand{\ben}{\begin{eqnarray}\displaystyle}
\newcommand{\een}{\end{eqnarray}}

\newcommand{\sectiono}[1]{\section{#1}\setcounter{equation}{0}}

\newtheorem{thm}{Th\'eor\`eme} [section]
\newtheorem{lemma}[thm]{Lemme}

\newtheorem{propo}[thm] {Proposition}

\newtheorem{rem}[thm]{Remarque}

\begin{document}

\centerline{\Large \bf  Quelques remarques sur le  champ des courbes  
 }\vskip .1cm
\centerline{\Large \bf  hyperelliptiques lisses en caract\'eristique  deux}

\medskip

\vspace*{4.0ex}

\centerline{\large \rm
 Jos\'e Bertin }

\vspace*{4.0ex}

\centerline{ Institut Fourier}
\centerline{   Universit\'e de Grenoble 1, 100 rue des Maths, Saint-Martin d'h\`eres F}
\vspace*{2.0ex}
\centerline{Jose.Bertin@ujf-grenoble.fr}

\vspace*{5.0ex}

\centerline{\bf Abstract} \bigskip

  In this note we describe the stack ${\cal H}_g$ of smooth hyperelliptic curves of genus $g$ over an algebraically closed field of characteristic two, as a quotient stack of  a smooth variety of dimension $3g+5$ by a non reductive  algebraic group, extending a well known result of Vistoli. As an application, we show the Mumford-Vistoli description of the Picard group of the stack  ${\cal H}_g$ is valid in this characteristic. We also  describe the natural stratification  of ${\cal H}_g$ by means of  higher ramification data. We point out  that after stable compactification $\overline {\cal H}_g$ is not smooth in codimension at least two.

\vfill \eject

\baselineskip=16pt

\tableofcontents

\sectiono{Introduction}

   Dans cette note, on souhaite mettre en \Žvidence quelques propri\Žt\'es du champ des courbes hyperelliptiques  en caract\'eristique deux, naturellement  li\Žes \ˆ \/ la sp\Žcificit\Ž \/ de cette caract\Žristique.  Ces remarques sont  toutes en relation  avec le caract\re Artin-Schreier du quotient par l'involution hyperelliptique. Notons qu'en genre un et caract\Žristique deux,  les points d'ordre deux d'une courbe elliptique  se r\Žduisent  \ˆ \/ un ou deux, et que sous les m\^emes conditions, en genre $g = 2$,  les  points de Weierstrass se r\'eduisent \ˆ \/ $3, 2$ ou $1$.  Rappelons  que pour le genre deux,    Igusa \cite {i} a d\'ecrit  d'une part l'espace des modules grossiers ${\text M}_2$ sur $\mathbb Z$, et  d'autre part  mis en \'evidence les particularit\'es li\'ees \`a \/ la caract\'eristique deux (forme de Rosenhain g\'en\'eralis\'ee, automorphismes).  Par exemple, les courbes de genre deux se r\'epartissent   en trois familles, d\'ecrites par une \'equation non homog\`ene de la forme Artin-Schreier:
 \be  Y^2 - Y = \begin{cases} \alpha X + \beta X^{-1} + \gamma (X-1)^{-1} \;(\alpha,\beta,\gamma \ne 0)\\
 X^3 + \alpha X + \beta X^{-1} \;(\beta \ne 0) \\
 X^5 + \alpha X^3. \end{cases} \ee
 Le premier   r\'esultat  de cette note est  une description du champ ${\cal H}_g$ sur $\overline {\mathbb F}_2$ comme champ quotient,  r\'esultat  qui  est enrichi  par  le fait que  ${\cal H}_g$  est pour cette caract\'eristique naturellement stratifi\Ž \/  au moyen du conducteur de Hasse, c'est \ˆ \/ dire de la ramification sup\'erieure.  Si $g=2$, c'est la stratification  de Igusa. On notera cependant que     ${\cal H}_1 \ne {\cal M}_{1,1}$ (voir Remarque 4.3).  
 En caract\'eristique $p \geq 0, \,p\ne 2$, il est bien connu que   si $g=1$ \cite {mum} ,  ou $g=2$ \cite {vg}
 \be\label{pich} \Pic ({\cal H}_g)  = \begin{cases} \frac {\mathbb Z} {12\mathbb Z}, & \text {si $g=1$}\\
 \frac {\mathbb Z} {10\mathbb Z}, &\text {si $g=2$}.\\ \end{cases}\ee
On notera  aussi que si $g=1$,  et  bien que les champs soient distincts,  on a l'\'egalit\Ž \/ $\Pic ({\cal H}_1) = \Pic ({\cal M}_{1,1})$ (Remarque 4.3).  La d\'etermination de $\Pic ({\cal H}_g)$ a \Ž t\Ž \/ \Ž tendue \ˆ \/ un genre arbitraire, et caract\'eristique $p\ne 2$,  par Arsie et Vistoli \cite {av};  de mani\`ere  un peu plus g\'en\'erale ils consid\`erent des rev\^etements \ˆ \/ groupe de Galois cyclique et diagonalisable (appel\'es par ces auteurs rev\^etements cycliques simples). 

Le  second r\'esultat  de cette note est la description de   $\Pic ({\cal H}_g)$ si $g \geq 1$, et toujours sur $\overline {\mathbb F}_2$.  Il d\'ecoule  facilement de la description du champ ${\cal H}_g$ comme champ quotient $[X/G]$, o\ \/ $X$ est un ouvert d'un espace de repr\'esentation de $G$, comme cela a \'et\Ž \/  remarqu\Ž \/ par  Vistoli \cite {av},\cite {vg}, .  La diff\'erence  dans notre cas venant  principalement du fait que  pour le groupe de Galois   $  \frac {\mathbb Z} {2\mathbb Z} \ne \mu_2$, les rev\^etements sont maintenant du type Artin-Schreier au lieu du type Kummer.  La description  obtenue peut \^etre vue comme  une variation sur le th\`eme des coefficients $a_0,a_1,\ldots,a_6$ de la forme de Weierstrass g\'en\'eralis\'ee en caract\'eristique deux \cite {si}. On remarquera  toutefois que  ${\cal H}_g$ est  encore un champ lisse, fait exceptionnel, non partag\Ž \/ en g\'en\'eral par les champs de rev\^etements en caract\'eristique positive \cite {bm}, \cite {rom}.  Il est remarquable que l'on puisse de cette mani\`ere retrouver ce fait qui  rel\`eve  essentiellement  d'arguments de d\'eformations \'equivariantes.  On remarquera cependant que la compactification stable $\overline {\cal H}_g$ bien que lisse en codimension un n'est pas lisse en codimension au moins deux si $g \geq 3$. 

On peut imaginer  qu'une description  analogue  en caract\'eristique  $p$  du champ de Hurwitz classifiant les rev\^etements $p$-cycliques   de $\mathbb P^1$   comme champ quotient, devrait  conduire \ˆ \/ des informations pr\' ecises sur  la nature de ses singularit\'es \cite {bm}.  Dans tout le texte, et sauf mention du contraire, un corps $k$ est alg\'ebriquement clos de caract\'eristique deux.
\medskip

\sectiono { Quelques calculs de Torseurs }
 \subsection{Objets de marque locale donn\'ee}
Dans cette section on rassemble des r\'esultats pour la plupart  bien connus, mais dans une formulation commode pour la suite.  De mani\re g\Žn\Žrale, on sait  que la classification des objets au dessus de la cat\Žgorie des $k$-sch\'emas  ($k$ corps, ou anneau de base), localement fpqc d'une  marque \emph {donn\Že} (confer. Demazure-Gabriel \cite {dg}, Ch III, ¤5),  est \' equivalente \ˆ \/  la donn\'ee du champ classifiant {\bf B}G, avec  pour groupe structural  G  le groupe des automorphismes de l'objet marque.  

Une d\Žfinition  pr\Žcise et plus g\Žn\Žrale est comme suit.  Soit $\cal M$ un champ alg\'ebrique  au dessus de $k$, i.e. au-dessus de la cat\Žgorie des $k$-sch\'emas. 
Soit $X$ un $k$-sch\Žma, et soit $P\in {\cal M}(X)$ un objet distingu\Ž \/, l'objet \emph {marque}. Disons que  $Q \in {\cal M}(U)$ est de {\it marque} $P$, s'il y a un morphisme $\alpha: U\to X$, un recouvrement (\Žtale) $U'\to U$, tels que 
\be Q\times_U U' \cong P\times_X U',\quad  \text {donc si }\quad \Isom (Q , P\times_X U) (U') \ne \emptyset\ee

On d\Žfinit un champ ${\cal M}_P$ (le champ des objets de marque $P$), dont les objets au-dessus de $U$ sont les couples $(Q,\alpha: U\to X)$, comme d\'efinis au-dessus, cette d\'efinition traduisant le fait que $Q$ est de marque $P$. Un morphisme $(Q,\alpha) \to (Q',\alpha')$ est un morphisme $f: Q\to Q'$ au-dessus de $u: U\to U'$, avec $\alpha' u = \alpha$. Ce champ n'est en fait comme on va le voir qu'une  pr\'esentation diff\'erente du champ classifiant ${\bf B} (G/X)$, o\ \/  $G = \Aut_X (P)$ est le $X$-groupe alg\Žbrique des automorphismes de $P$, suppos\Ž \/ lisse \cite {lmb}.
\begin{propo}
  Sous les hypoth\ses qui pr\Žc\dent,  on a un isomorphisme de champs
$${\cal M}_P \stackrel\sim \longrightarrow ( {\bf G}/X)$$
$G = \Aut _k(P)$ \'etant comme indiqu\Ž,  le  $X$-groupe alg\'ebrique des automorphismes de $P$.  
\end{propo}
\begin {proof} Rappelons bri\vement la d\Žfinition du morphisme ${\cal M}_P \stackrel\sim  \rightarrow  {\bf G}/X)$ (voir \cite {dg} pour des d'etails). Si $(Q\in {\cal M}(U), \, f:U\to X)$ est un objet de marque $P$, il est clair que le $U$-sch\Ž ma $\Isom_U (Q , P\times_X U)$ est  de mani\'ere naturelle, c'est \ˆ \/ dire pour l'action \ˆ \/ droite de $G\times_X U$,  un  $ G$ fibr\Ž \/ principal de base $U$, i.e. un $G\times_X U$-torseur. Un morphisme entre objets de marque $P$ induit de mani\`ere \'evidente un morphisme  entre torseurs associ\'es.  En sens inverse, si on a un morphisme $\alpha: U\to X$, et si $Q \to U$ est un $G\times_X U$-torseur,  le produit contract\Ž \/  \footnote { 
 Supposons le torseur $Q\to U$, d\'efini par un cocycle $g_{i,j}: U_{i,j} \to G\times_X U$, relativement \ˆ \/ un recouvrement \'etale $(U_i\to U)_i$ de $U$. On peut alors voir $(g_{i,j})$ comme une donn\'ee de descente sur la collection $P_i = P\times_X U_i$, relativement \ˆ \/ ce recouvrement. Cela d\'efinit un objet  $x\in {\cal M}(U)$ tel que $Q = \Isom (x , P\times_X U)$.}
\be Q \times^{G\times_X U}  \;(P\times_X U)\ee
d\Žfinit  un objet de marque $P$ de base $U$; de mani\re plus pr\Žcise $(P,\alpha) \in {\cal M}_P(U)$. Ces deux constructions sont clairement inverses l'une de l'autre. D'une autre mani\`ere on peut observer que  le champ ${\cal M}(P)$ est une gerbe  sur $X$. cette gerbe est  neutre car $P$ en est une section sur $X$. Le r\'esultat d\'ecoule alors de (\cite {lmb}, Lemme (3.31))\end{proof}

\subsection{ Le champ classifiant  $\Bb G$}
 
Un exemple  simple  qui illustre la construction pr\Žc\Ždente revient \ˆ \/  prendre comme marque $k^n$,  de sorte que  $\cal M$  est dans ce cas le champ des faisceaux localement libres de rang $n$ sur les $k$-sch\Žmas.  C'est le champ alg\Žbrique ${\cal F}ib_{n,X/S}$ de  (\cite {lmb}, Th\Žor\me (4.6.2.1)), avec $X = S = \Spec k$.  Les sections au dessus de $U$ sont les fibr\'es vectoriels de rang $n$ sur $U$, les morphismes \'etant les morphismes de changement de base. Alors le principe qui vient d'\^etre rappel\Ž \/   dit  que ce champ est  isomorphe \ˆ \/ $\Bb \GL (n)$. Plus significatif pour la suite est le cas du groupe projectif lin\'eaire $\PGL (n) = \Aut (\mathbb P^n)$. Le champ classifiant $\Bb \PGL (n)$ est   isomorphe au champ dont les objets sont les fibrations en $P^n$ au dessus d'un $k$-sch\'ema (sch\'emas de Severi-Brauer). 

Le cas qui nous sera  utile pour la suite est un m\Žlange de ces  deux exemples.  Fixons un faisceau localement libre  $V$   sur $P^1,$ de rang $r \geq 1,$ donc en vertu d'un th\Žor\me bien connu de Grothendieck, de la forme
\be V = {\cal O}(n_1)\oplus \cdots \oplus {\cal O}(n_r)\ee
 avec $n_1\leq \cdots \leq n_r$. Soit  la cat\Žgorie fibr\Že en groupo\"{\i}des  $\cal P$, dont les objets sont les couples $(D\to S, E)$, avec $D\to S$ une fibration en $\mathbb P^1$, et $E$ un fibr\Ž \/ vectoriel sur $D$ localement isomorphe \ˆ \/ $V$,   c'est \ˆ \/ dire \Žtale-localement sur $S$ de la forme $\oplus_{i=1}^r {\cal L}_i$, o\ \/ ${\cal L}_i$ est inversible de degr\Ž \/ relatif $n_i$.  Noter que cette d\Žfinition  a un sens car localement pour la topologie \Žtale $D = \mathbb P_S^1$. L'identification des objets dans cet exemple est un exercice ais\Ž \/\footnote {On peut moduler la d\Žfinition, par exemple en ajoutant   \ˆ \/ $V$ une section globale. Les objets sont alors des triples $(D ,E,s), \, s\in \Gamma (D,E)$. Les automorphismes sont dans ce cas assujettis \ˆ \/ fixer la section.}:
 
 \begin{propo}
 \begin{enumerate}
\item La cat\Žgorie fibr\Že en groupo\"\i des $\cal P$ est un champ alg\Žbrique, i.e. ${\cal P} \cong \Bb G$. 
\item La marque   des objets  de ${\cal P}$ est $(\mathbb P^1,V)$. Le groupe  $G$ des automorphismes du mod\le est   le groupe qui lin\Žarise universellement le fibr\Ž \/  $V$.  Ce groupe  s'ins\re  dans une suite exact
\be 1 \rightarrow H = \Aut (V) \rightarrow G \rightarrow \PGL (1) \rightarrow 1  \ee
\end{enumerate}\end{propo}

\begin {proof} Le principe g\Žn\Žral  indiqu\Ž \/ au dessus s'applique si on plonge les objets dans un champ ambiant. Une preuve directe de l'identification 1) est aussi rapide. Le premier point se r\Žsume essentiellement \ˆ \/ montrer que si $(D,E)$ et $(D',E')$ sont deux objets de $\cal P$ au-dessus du sch\Žma $S$, alors le faisceau sur $S_{et}$
\be  (T\to S) \mapsto \Isom_T \left ((D_T,E_T) , (D'_T,E'_T)\right )\ee
est repr\Žsentable  par un $S$-sch\Žma s\Žpar\Ž \/ de type fini.  On peut proc\Žder en deux \Žtapes. Soit d'abord le foncteur $\Isom_S (D,D')$ qui est  repr\'esent\Ž \/ par  un $S$ sch\'ema  affine. Si $D = \mathbb P^1_S$, c'est un torseur $\cal I$ sous le groupe $\PGL (1)_S$. Si $\Phi: D\times _S {\cal I} \buildrel\sim\over\rightarrow D'\times_S \cal I$ est l'isomorphisme universel, on est ramen\Ž \/  \ˆ \/ prouver que le foncteur  d\Žfini sur la cat\Žgorie des $\cal I$-sch\'emas 
 
\be (U\to {\cal I}) \mapsto \Isom_U (E_U, \Phi^*(E')_U)\ee
 
est repr\'esentable.  L'argument est bien connu (voir par exemple \cite {lmb}, Th\'eor\ me 4.6.2.1);  ce foncteur 
est en fait repr\'esent\Ž \/  par un  sch\'ema affine sur $\cal I$. Comme il est manifeste que l'objet $(\mathbb P^1,V)\in {\cal P}(k)$ d\'efinit un atlas,   le premier point est essentiellement clair. 

 Identifions  maintenant  la marque des objets de $\cal P$. Il est clair que toute fibration en $\mathbb P^1$ de base $S$ est localement (pour la topologie \Žtale) de la forme $\mathbb P^1_S$, de mani\re \Žquivalente, a localement pour la topologie \'etale une section. Supposons donc $D = \mathbb P^1_S$. Dans ce cas  le  module localement libre  $E$ \Žtant   localement isomorphe \ˆ \/  $V$, on peut  quitte \ˆ \/ localiser si n\Žcessaire, supposer  que 
 $$E =  {\cal O}_{\mathbb P^1_S}(n_1)\oplus \cdots \oplus {\cal O}_{\mathbb P^1_S} (n_r) = p^* (V)$$
 en d\'esigant par  $p: \mathbb P^1_S \to \mathbb P^1_k$  le morphisme de changement de base.  Ce qui  montre que  la marque est bien le couple $(\mathbb P^1,V)$.

Pour d\Žcrire le groupe des automorphismes de l'objet marque,  rappelons que si un groupe alg\Žbrique $G$, de multiplication $\mu$ agit sur une vari\Žt\Ž \/ $X$, l'action \Žtant not\Že $\sigma: G\times X \to X$, une $G$-lin\Žarisation sur un faisceau coh\Žrent $F$, est un isomorphisme
\be \Phi:  \sigma^* (F) \buildrel\sim\over\longrightarrow p_2^* (F) \ee 
v\Žrifiant la relation de cocycle  
\be (\mu,1)^* \Phi = p_{2,3}^* \Phi . (1,\sigma)^*\Phi \ee
 En g\Žn\Žral un faisceau  $F$ localement libre de rang $n$, invariant par $G$, donc tel que $g^*(F)\cong F$ pour tout $g\in G$,  n'est pas $G$-lin\Žarisable, cependant il le devient si le groupe $G$ est agrandi convenablement.  De mani\re plus pr\Žcise, l'argument rappel\Ž \/ en d\Žbut de preuve montre qu'il existe un groupe alg\Žbrique $\tilde G$, et un morphisme de groupes
$\tilde G \to G$, tel que les $G$-lin\Žarisations de $F$ correspondent de mani\re bijective aux sections de $\tilde G \to G$.  En fait il suffit d'appliquer  la construction  du d\Žbut au couple de faisceaux $(\mu^*(F), p_2^*(F))$,   $\mu$ et $p_2$ \Žtant respectivement  l'action,  et la projection $G\times X \to X$, de sorte que: 
\be \tilde G = \Isom (p_2^* (F) , \mu^* (F))  \ee
 D'une mani\re simplifi\Že, les points de $\tilde G$ sont les couples $(g,\phi)$, avec $g\in G$, et $\phi$ un isomorphisme $ F\buildrel\sim\over\rightarrow g^*(F)$. La loi de composition  est $(g,\phi)(h,\psi) = (gh,h^*(\phi)\psi)$. Si on a pour tout $g\in G, \;g^*(F) \cong F$, alors  ce groupe s'ins\re dans une extension, g\Žn\Žralisant 1.1:
\be1\rightarrow \Aut (F) \rightarrow \tilde G \rightarrow G\rightarrow 1 \ee
 Revenons \ˆ \/ la situation de d\Žpart, donc $F = V, \, X = \mathbb P^1$. Soit $\tilde Q \to S$ un $\tilde G_S$-torseur. Avec ce torseur on construit en premier un fibr\Ž \/ en $\mathbb P_1$, qui est le fibr\Ž \/ associ\Ž \/ 
 $$P ={ \tilde P \times_S}^{\tilde G} \, \mathbb P^1 \rightarrow S $$
 On construit ensuite un fibr\Ž \/ vectoriel $E\to P$, localement sur $S$ isomorphe \ˆ \/ $V$, en utilisant le fait que par construction $\tilde G$ agit sur le fibr\Ž \/ $V \to \mathbb P^1$. On prend le fibr\Ž \/ vectoriel associ\Ž \/:
 $$ E = {\tilde Q \times_S}^{\tilde G} \, V_S \to P  $$
 Cette construction jointe \ˆ \/ la pr\Žc\Ždente conduit \ˆ \/ l'identification $\cal P  \cong \Bb \tilde G$, et donc donne le r\Žsultat en ajustant les notations.
 \end{proof}

Notons que l'extension 2.10  n'est  en g\Žn\'eral pas scind\Že,  du fait que ${\cal O}(n)$ n'est $\PGL (1)$-lin\Žarisable  que  si  $n$ est pair; il poss\de  bien entendu une $\GL (2)$-lin\Žarisation canonique. En effet soit $\alpha_g: {\cal O}(1) \buildrel\sim\over\rightarrow g^* ({\cal O} (1))$ la lin\Žarisation tautologique de ${\cal O} (1)$,  la lin\Žarisation canonique de ${\cal O}(n)$ relative \ˆ \/ $\GL(2)$  est celle donn\Že par $\alpha_g^{\otimes n}$. Cela d\Žfinit la lin\Žarisation tautologique de $V$:
\be  \oplus_{i=1}^r \,  \alpha_g^{\otimes n_i}:  V \buildrel\sim\over\rightarrow g^* (V) \quad (g\in \GL (2)) \ee
Noter que pour cette action $\lambda 1_2$ agit diagonalement  sur $V$ par $ \text {diag} (\lambda^{n_1}, \cdots ,\lambda ^{n_r}). $
 
 On notera cependant  que la lin\Žarisation tautologique  de $V$ sous $\GL (2)$ se descend dans tous les cas  en une  lin\Žarisation sous $\GL (2) / \mu_d$, en notant $d $  le pgcd des $n_i$.   
 
 Le groupe  $H$ de la proposition 1.2 est ais\Ž \/ \ˆ \/ d\'ecrire.  Supposons la partition $(n_i)$ de $r=\sum_{i=1}^r n_i$,  prenant $p$ valeurs distinctes $k_1 < \ldots < k_p$,  $k_i$ apparaissant $\alpha_i$ fois.  Alors $V$ a une filtration canonique $F^\bullet$ (de Harder-Narasimhan)  
 $$0 \ne F^p \subset \cdots \subset F^1 = E,\quad  F^j/F^{j+1} = {\cal O}_{\mathbb P^1}(k_j)^{\oplus  \alpha_j}  $$
 En fait si $V = \oplus_{j=1}^p{ \cal O}_{\mathbb P^1}(k_j)^{\oplus \alpha_j} $, on a 
 $F_j =  \oplus_{m=j}^p {\cal O}_{\mathbb P^1}(k_m) ^{\oplus (\alpha_m)}. $ 
En cons\'equence  un automorphisme de $V$ respecte  la filtration $F^\bullet$.   Il est  alors imm\Ždiat  que $H = \Aut (V)$ s'identifie \ˆ \/  un  groupe produit semi-direct 
$$ H = U \rtimes \prod_{j=1}^p \GL (\alpha_i)  $$
expression dans la quelle le groupe unipotent $U$, est un sous-groupe du groupe   des automorphismes de $V$ induisant l'identit\Ž \/ dans le gradu\Ž \/ associ\Ž \/ $gr^\bullet (F^\bullet)$. 
 
 \begin{rem}\end{rem}
Sous certaines conditions $G$ est aussi un produit semi-direct. Supposons par exemple $\pgcd (k_1\alpha_1,\cdots ,k_p\alpha_p) = d$. Soit alors  une relation de Bezout $d = \sum_j m_j   k_j\alpha_j  $. Consid\Žrons  le caract\re  de $H$ donn\Ž \/ par
$$\chi :  H \rightarrow \prod_{j=1}^p \GL (\alpha_i)\buildrel\psi\over\rightarrow G_m  $$
avec $\psi ((g_j)) = \prod_j  \det (g_j)^{m_j}$, et soit $K $ le noyau de $\chi$.   
 Sous ces conditions  le groupe $H$   est  le  produit semi-direct  $H =  K\rtimes \GL (2) /\mu_d$.     Notons $\alpha: \GL (2) \to G$ la section induite par la lin\Žarisation canonique sur chaque facteur. Alors  les \Žl\Žments de $K$ qui sont  dans l'image de $\alpha$, sont  repr\Žsent\Žs par les matrices diagonales $\text {diag}(\lambda^{n_1}, \cdots ,\lambda^{n_r})$, avec $\psi (\alpha (g)) = \lambda^{\sum m_in_i} = \lambda^d = 1$.  Cette intersection est donc r\Žduite \ˆ \/ l'identit\Ž \/.  
 
 Dans le cas g\Žn\Žral, et $d$ \Žtant  suppos\Ž \/ pair, la suite (1.1) est scind\Že.  
 Notons en effet que si les $n_i$ sont tous pairs   on obtient une section de $G \to \PGL (1)$ par factorisation de $g\in \GL (2) \mapsto \oplus_i\,  \alpha_g^{\otimes n_i} \det (g)^{-n_i/2}$.  $\lozenge$
 
  \sectiono   {Courbes hyperelliptiques  en caract\'eristique deux}
  \subsection {Le champ ${\cal H}_g$}
  
   On fixe  dor\'enavant un corps  $k$ alg\'ebriquement clos de caract\'eristique deux, par exemple $k = \overline {\mathbb F}_2$.  Les sch\'emas sont des $k$-sch\Žmas  de type fini. Rappelons  qu'une courbe hyperelliptique $p: C \to S$,  de genre $g \geq 2$, de base  un sch\Žma $S$,  est une $S$-courbe projective lisse \ˆ \/ fibres connexes, munie d'un $S$-automorphisme involutif $\tau: C \stackrel\sim \rightarrow  C$, de sorte  que toute fibre g\'eom\'etrique  $C_s$  est une courbe hyperelliptique de genre $g$,  d'involution hyperelliptique la restriction de    $\tau$.   L'involution $\tau$, qui est alors unique, est appel\'ee l'involution hyperelliptique de la courbe $p: C\to S$.    La courbe quotient $D = C/\tau$ est  une  fibration  tordue en $\mathbb P^1$ de base $S$, voulant dire que les fibres g\'eom\'etriques sont des $\mathbb P^1$). Notons que  le quotient $\pi: C \to D$ est plat de degr\Ž \/ deux.  Il est connu que la formation de $D$ commute \ˆ \/ tout changement de base, du fait que les courbes sont lisses (voir par exemple \cite {kl}). Justifions d'abord que les courbes hyperelliptiques (lisses) de genre $g\geq 1$ forment un champ alg\'ebrique.  Si  on accepte le cas $g=1$, on notera que l'involution n'est plus unique, cependant  la d\'efinition  du champ garde un sens. Une courbe (hyper)elliptique de genre un est donc un couple $(C,\tau)$.   
Soient $p: C\to S$ et $p': C' \to S'$, deux courbes hyperelliptiques de base $S$. Un morphisme au dessus de $h: S \to S'$ est comme d'habitude d\'efini par un diagramme cart\'esien
  \be  \CD
     C     @>  {f} >>  C' \\ @V pVV @VV {p'}V \\
     S  @>h>>  S' \endCD  \ee
le morphisme $f$  v\'erifiant alors de mani\'ere automatique la r\'egle de commutation
    $\tau' f = f\tau$. Les courbes hyperelliptiques forment  une cat\'egorie fibr\'ee en  groupo\"{\i}des   ${\cal H}_g$ ( ou ${\cal H}_{g,k}$) au dessus de la cat\'egorie des $k$-sch\'emas.  De mani\'ere plus pr\'ecise on a le r\'esultat  connu  et ais\'e suivant  \cite {bm},\cite {rom} que nous \'enoncons sans preuve:  
    \begin{propo}  Pour tout $g\geq 1$, le groupo\"{\i}de ${\cal H}_g$ est un champ alg\'ebrique  de Deligne-Mumford  lisse. \end{propo} 
    Une remarque cependant s'impose. La lissit\Ž \/ du champ ${\cal H}_g$ est  un fait exceptionnel, qui n'est pas partag\Ž \/ par beaucoup d'autres champs de Hurwitz en caract\'eristique positive,  du moins lorsque la caract\'eristique  de $k$ divise l'ordre du groupe de monodromie. Ce fait  d\'ecoule par exemple du  r\'esultat plus pr\'ecis  qui dit que  l'anneau versel des d\'eformations $\frac {\mathbb Z} {2\mathbb Z} $-\'equivariantes d'un point fixe formel (i.e. $k[[t]]$) de conducteur de Hasse $m$  ($m$ est forc\'ement impair) est lisse de dimension $\frac { m+1} 2$  \cite {bm}.

    Nous allons  proposer une description beaucoup plus  explicite  de ${\cal H}_g$,   et se pr\^etant mieux   au ''calcul'', description qui rendra par ailleurs la lissit\Ž \/ \'evidente. Soit $p: C \to S$ un objet de ${\cal H}_g$, et  soit $\pi: C\to D$ le rev\ tement associ\Ž \/ de degr\Ž \/ deux.  La  $S$-courbe $q: D \to S$, est  un fibr\Ž \/ en coniques, donc localement  pour la topologie fpqc (en fait \'etale localement)   une droite projective $\mathbb P^1_S$.  Soit par  ailleurs le faisceau localement libre de rang deux ${\cal E} = \pi_\star ({\cal O}_C)$. Il est \Žquip\Ž \/ d'une part  d'une involution naturelle, celle induite par $\tau$, not\'ee  pour cette raison $\tau$, et aussi d'une structure de $\cal  O_D$-alg\bre; de sorte qu'on r\'ecup\`ere $C$ par l'op\'eration $C = \Spec _D ({\cal E})$.  
    
    Soit la suite exacte  d\'efinie en faisant le quotient par ${\cal O}_D = \ker (\tau - 1)$:
    \be 0\rightarrow {\cal O}_D \rightarrow \cal E \rightarrow \Cal L \rightarrow 0  \ee
    Soient  $m_i \, (i=1,\cdots,b)$   les conducteurs de Hasse locaux en les $b$ points de branchement de courbe hyperelliptique $\pi: C\to D = \mathbb P^1$ de genre $g$,  d\Ž finie sur le corps $k$.  On notera que la formule de Riemann-Hurwitz donne 
    \be g+1 = \sum_{i=1}^b  \frac {m_i+1} 2\ee
     Dans la suite on posera $m = g+1$.
    \begin {lemma} Soit une courbe hyperelliptique de base $S$. Localement sur $S$  on a $D \cong \mathbb P^1_S, \,\,{\cal L} \cong {\cal O}_{\mathbb P^1_S} (-m)$, de plus la suite exacte (2.2) est (localement sur $S$) scind\'ee. \end{lemma}    
    \begin{proof}   On a par le th\'eor\`eme de Riemann-Roch 
    $\chi ({\cal O}_C) = \chi ({\cal E}) = \deg ({\cal L}) + 2$,
     de sorte que $\deg ({\cal L}) = -m$.  Pour montrer que (2.2) est localement scind\Že \/ sur $S$,  on peut supposer $S$ affine,  $D = \mathbb P^1_S$, et que ${\cal L} = {\cal O}(-m)$. Alors 
     $$\Ext^1 ({\cal  O}(-m), {\cal O}_D) = H^1(D,{\cal O}(m)) = H^1(S,p_*{\cal O}(m)) = 0$$
      Le r\'esultat en d\'ecoule. \end{proof}
    
    Dans la suite,  lorsque les conclusions du lemme sont satisfaites, donc sur un recouvrement fpqc de $S$, on dira que la courbe hyperelliptique est rigidifi\'ee\footnote { En caract\'eristique $p\ne 2$, l'extension (3.2) a une section canonique, du fait de la d\'ecomposition de $E$ en sous-espaces propres de $\tau$, mais cela n'est plus le cas si $p = 2$.} d\`es lors que  des identifications $D\cong \mathbb P^1_S$ et ${\cal E} \cong {\cal O}_D \oplus {\cal O}_D (-m)$ (i.e. une section de (2.2)) sont choisies. Dans ce cas la courbe $p: C\to S$ est  d\'etermin\'ee de mani\`ere unique par la structure de ${\cal O}_D$-alg\bre sur ${\cal E} = \pi_* ({\cal O}_C)$, c'est \ˆ \/ dire par la donn\'ee de deux sections $A\in \Gamma ({\cal O}_D (m)$ et $B\in \Gamma ({\cal O}_D (2m)$. 
    
 De mani\ re  plus pr\Ž cise, pour toute base locale $e$ de ${\cal L} = {\cal O}_D (-m)$, le carr\Ž \/ de $e$ est
    \be  e^2 = B(e^{\otimes 2}) \,+ \, A(e) e \ee
Pour que la structure d'alg\`ebre (2.3) d\'efinisse en tout point $s\in S$ une courbe hyperelliptique, il est n\'ecessaire que  la condition $A(s) \ne 0$ soit satisfaite. Dans le cas contraire la projection sur $\mathbb P^1$ est radicielle; de m\^eme si $B(s) = 0$, la fibre est r\'eductible.  Les conditions sur le couple $(A,B)$ qui assurent que la courbe $C(A,B) = \Spec_D ({\cal E}(A,B))$ est hyperelliptique lisse seront explicit\'ees ci-dessous. Dans cette notation, on d\'esigne par ${\cal E}(A,B)$ l'alg\bre  d\'efinie par la multiplication (2.3) sur ${\cal O}_D \oplus {\cal O}_D(-m)$. Dans ce cas  on notera que l'automorphisme hyperelliptique $\tau$ de $C(A,B)$  a pour expression \/ par $\tau (e) =  A(e) + e$. 

On notera dans la suite $R\subset D$ le diviseur de ramification (relatif) du rev\^etement $\pi: C \to D$. Notons toujours sous les m\^emes conditions  le fait \'el\'ementaire:
    \begin{lemma} \begin{enumerate}
     \item Avec les notations introduites au-dessus, le  diviseur $F$ des points fixes de $\tau$ (ou diviseur de Weierstrass)  a pour \'equation $A = 0$.
 \item     On a l'\'egalit\Ž  \/ $F = R$, et  ${\cal O}_C(R) \cong \pi^* ({\cal L}^{-1})$. \end{enumerate}\end{lemma}
    \begin{proof}   1)  Localement, le diviseur des points  fixes de $\tau$ a un id\'eal  engendr\Ž \/ par
     les \'el\'ements $\tau (\xi) - \xi$, pour $\xi$ section de  ${\cal O}_C$. Il est clair que cet id\'eal est engendr\Ž \/ par $\tau (e) - e =  A(e)$. \\
     2) Noter que pour un rev\^etement galoisien de groupe $G$ entre courbes lisses, et dans le contexte de la ramification sauvage, le diviseur de ramification, c'est \ˆ \/ dire  d'un point de vue local  l'exposant de la diff\'erente, s'exprime en fonction  de la fonction d'ordre, donc des multiplicit\'es des diviseurs de points fixes $F_\sigma$ des \'el\'ements  $\sigma \in G$ par la formule bien connue \cite {rom},\cite {jps}:
     \be R = \sum_{\sigma \in G,\sigma \ne 1} \, F_\sigma\ee
Dans notre situation $G $ est d'ordre deux, cela donne l'\'egalit\Ž \/ $R = F$.  Le dernier point d\'ecoule imm\'ediatement de la description locale (1).   \end{proof}
   \subsection{Forme de Rosenhaim g\'en\'eralis\'ee} 
     
 Nous allons analyser les contraintes impos\'ees  sur le couple $(A,B) $ pour que $C(A,B)$ soit hyperelliptique  lisse.  Soit toujours  $S = \Spec (k)$,  avec $k$ corps alg\'ebriquement clos de caract\'eristique deux. On a $D = \mathbb P^1$ de coordonn\`ees homog\`enes  $X,Y$, et la courbe $C(A,B)$ est d\'ecrite sur les ouverts affines $Y\ne 0$ et $X\ne 0$, par les \'equations respectives\footnote { En genre $g=2$, Igusa \cite {i} utilise le mod\`ele birationnel  (forme de Rosenhaim) $xy^2 + (1+ax+bx^2)y + x^2(c+dx+x^2) = 0$ valable en toute caract\'eristique.}
\be z^2 + A(x,1)z = B(x,1)\quad \text {resp. } w^2 + A(1,y)w = B(1,y) \quad (w = \left(\frac y x\right)^m z)  \ee

avec 
$$A(X,Y) = a_0X^m+ \cdots +a_mY^m, \quad  B(X,Y) = b_0X^{2m} + \cdots + b_{2m}Y^{2m}$$
 formes binaires, respectivement de degr\Ž \/ $m$ et $2m$.
La lissit\Ž \/ est   visiblement \'equivalente au fait que les deux syst\`emes, dans lesquels l'indice $X$ (resp. $Y$) d\'esigne une d\'eriv\'ee partielle 
\be  \begin{cases} A(x,1) =   A'_X (x,1) ^{2} B(x,1) + B'_X  (x,1)^2 = 0\\   A(1,y) =   A'_Y (1,y) ^{2} B(1,y) + B'_Y  (1,y)^2 = 0\end{cases}  \ee
 n'ont pas de solutions $x$ (resp. $y$).  De mani\`ere homog\`ene,  cela veut dire que les syst\`emes
\be \begin{cases} A(X,Y) =   A'_X (X,Y) ^{2} B(X,Y) + B'_X (X,Y)^2 = 0 \\  A(X,Y) =   A'_Y (X,Y) ^{2} B(X,Y) + B'_Y (X,Y)^2 = 0 \end{cases} \ee
n'ont pas de solutions $(x,y)$ avec $y\ne 0$ pour le premier  (resp. $(x,y), \;x\ne 0$ pour le second). Dans la suite, on identifiera l'espace affine des couples $(A,B)$ \ˆ \/ $\mathbb A^{m+1}\times \mathbb A^{2m+1} = \mathbb A^{3m+2}$. Les coordonn\'ees \'etant $(a_0,\cdots ,b_{2m})$. On consid\`erera   l'espace affine $\mathbb A^{3m+2}$ comme muni de plusieurs actions naturelles,  d'une part de l'action de $\mathbb G_m^2$
$$ (t,s).(A,B) = (tA,sB), \quad (t,s)\in ( k^*)^2 $$
d'autre part de l'action du groupe unipotent $\mathbb G_a^{m+1}  = \Gamma (\mathbb P^1 , \Cal O (m))$ (groupe additif)\footnote { L'action de $\mathbb G_a^{m+1}$ n'est pas normalis\'ee par $\mathbb G_m^2$, mais seulement par le sous-groupe \ˆ \/ $1$-param\`etre $\lambda \mapsto (\lambda,\lambda^2)$. Le groupe produit semi-direct correspondant a une interpretation claire, comme groupe assurant des changements de rep\`eres (voir \S 4).}
\be \alpha .(A,B) = (A,B+\alpha A+\alpha^2)  \ee

Noter que cette action ne fait que traduire  un changement de section dans  la suite (3.2). Enfin l'action de $\GL (2)$ qui traduit un changement lin\'eaire simultan\Ž \/ des coordonn\'ees
   $(A,B) \mapsto (A\circ \sigma^{-1} , B\circ \sigma^{-1}). $ 
 
Pour transcrire en des termes commodes la condition de lissit\Ž \/, notons d'abord le lemme \'el\'ementaire suivant, dans lequel on note   $\Res  = \Res _{m,4m-2}$ le r\'esultant de deux formes binaires de degr\'es respectifs $m$ et $4m-2$. Ce lemme montre que le  polyn\^ome $\Delta (A,B)$ jouit de propri\'et\'es analogues \ˆ \/ celles du discriminant \cite {gkz}.
\begin {lemma}\begin{enumerate} \item Il existe un polyn\^ome $\Delta (A,B)$ en les coefficients de $A$ et $B$, tel que   
\be \Res (A , {A'_X}^2 B  + {B'_X} ^2) = a_0 ^2 \; \Delta, \; \Res (A , {A'_Y} ^2 B + {B'_Y} ^2) = a_m ^{2} \;\Delta  \ee
Les coefficients $a_i$ de $A$ \'etant de poids un, et ceux de $B$, les $b_j$, \'etant de poids deux,  $\Delta (A,B)$ est homog\`ene de degr\Ž \/ $8m-4$.
\item Si $m$ est pair, on a 
\be \Delta = \Res (A , \frac {{A'_X} ^2 B  + {B'_X}^2}  {Y^2}) \;=\; \Res (A , \frac {{A'_Y}^2 B + {B'_Y}^2} {X^2})  \ee
\item $\Delta $ est  invariant par l'action de $\mathbb G_a^{m+1}$; par ailleurs une substitution lin\'eaire $\sigma \in GL(2)$ des variables donne $\Delta (A\circ \sigma , B\circ \sigma) = \det (\sigma)^{m(4m-2)} \Delta (A,B)$. \end{enumerate} \end{lemma}

\begin{proof}  1-  Utilisant le fait que $k$ est de caract\' eristique deux, et que $B$ est de degr\Ž \/ $2m$,  la formule d'Euler donne
$$  XA'_X + YA'_Y = mA, \,\, XB'_X = YB'_Y $$
ce qui conduit \ˆ \/ l'identit\Ž \/
$$X^2 \left ( {A'_X} ^2 B  + {B'_X}^2 \right ) = Y^2 \left ({A'_Y}^2 B + {B'_Y}^2 \right ) + m^2 A^2B  $$
Si on forme le r\' esultant de  $A$ avec le polyn\^{o}me homog\`ene  de degr\Ž \/ $4m$,  $X^2({A'_X} ^2 B  + {B'_X}^2),$ l'\'egalit\Ž \/ (3.12) donne  imm\'ediatement le r\'esultat (2.8). \\
2- Si $m$ est pair,   (3.12) montre que $X^2$ divise ${A'_Y}^2 B + {B'_Y}^2 $ et $Y^2$ divise $ {A'_X} ^2 B  + {B'_X}^2$, de sorte que dans ce cas, on bien  (2.9). \\
3-  Pour v\'erifier l'invariance de $\Delta$ notons que si $A^* = A$, et  $B^* = B+\alpha A+\alpha^2$, alors 
$${({{A^*}'}_X)}^2 B^* +{ ({{B^*}'}_X)}^2 = {A'_X}^2 B + {B'_X}^2 + \alpha A({A'_X}^2) + {\alpha '_X}^2 A^2$$
de sorte que l'\'egalit\Ž 
$$\Res (A^*, {{A^*}'}_X^2 B^* + {{B^*}'}_X^2) = \Res (A , {A'_X}^2 B + {B'_X}^2)$$
  d\'ecoule des propri\'et\'es \'el\'ementaires du r\'esultant (\cite {gkz}, Chapter 12). Pour prouver que $\Delta$ est un semi-invariant  de $\GL (2)$ de poids $m(4m - 2)$, il suffit de v\'erifier ce  fait, donc 
  $$\Delta (A\circ \sigma,B\circ \sigma) = \det (\sigma)^{m(4m-2)}) \Delta (A,B)$$
   pour $\sigma = \begin{pmatrix} 0 & 1\\ 1 & 0\end{pmatrix} , \, \begin{pmatrix} 1 & 1 \\ 0&1\end{pmatrix}, \, \begin{pmatrix} \lambda & 0 \\ 0 & \mu\end{pmatrix} $.  Dans le premier cas, avec la notation  $F^* (X,Y) = F(Y,X)$, on a $(F^*)'_X = (F'_Y)^*$;  on peut   donc \'ecrire 
   $$a_m^2 \Delta (A^*,B^*) = \Res (A^* , {(A^*)'_X}^2 B^* + {(B^*)'_X}^2 = \Res (A^* , ({A'_Y}^2 B + {B'_Y}^2)^*) $$
   $$= \Res (A , {A'_Y}^2 B + {B'_Y}^2) = a_m^2 \Delta (A,B)$$
   D'o\ \/ le r\'esultat dans ce cas.  Le second cas est tout aussi imm\'ediat. Dans le cas trois, avec la notation $F^*(X,Y) = F(\lambda X , \mu Y)$, on a du fait de la propri\'et\Ž \/ de semi-invariance du r\'esultant:
   $$\Res (A^* , {(A^*)'_X}^2 B^* + {(B^*)'_X}^2)  = \Res (A^* , \lambda^2 ({A'_X}^2 B + {B'_X}^2) $$
   $$ = \lambda^{2m} (\lambda\mu)^{m(4m - 2)} \Res (A , {A'_X}^2B + {B'_X}^2)$$
   Le r\'esultat d\'ecoule alors de la d\'efinition de $\Delta$.
 \end{proof}
Par exemple, dans le cas  $m=1$,  on trouve 
$$\Delta (A,B) = a_0^2b_2 + a_1^2b_0 + a_0a_1b_1 + b_1^2$$
  Si $m=2$, donc $g=1$, et en utilisant les conventions d'\'ecriture usuelles \cite {si} 
$$A(X,Y) = a_1XY + a_3 Y^2, \;\; B(X,Y) = X^3Y + a_2X^2Y^2 + a_4XY^3 + a_6 Y^4$$
on trouve\footnote {L'expression compl\te est:  $\Delta (A,B) = a_2^4(a_1^4b_0^2 + b_1^2) + a_0a_1^4a_2^3b_1^2 + a_0^2a_1^4a_2^2b_2^2 + a_0^3a_1^4a_2b_3^2 + a_0^4(a_1^4b_4^2 + b_3^4) + a_1a_2^3(a_1^4b_0b_1 + a_1^2b_1^3) + a_1^2a_2^2(a_1^4b_0b_2 + a_1^2b_1^2b_2) + a_2(a_1^3 + a_0a_1a_2)(a_1^4b_0b_3 + a_1^2b_1^2b_3) + a_1^4(a_1^2b_0 + b_1^2)(a_1^2b_4 + b_3^2) + a_0a_1^5a_2^2b_1b_2 + a_0a_1^6a_2b_1b_3 + a_0a_1^2b_1(a_1^4b_4 + b_3^2)(a_1^3 + a_0a_1a_2) + a_0^2a_1^5a_2b_2b_3 + a_0^2a_1^4(a_1^2b_4 + b_3^2) + a_0^3a_1^3b_3(a_1^2b_4 + b_3^2). $}
 pour $\Delta$ l'expression bien connue:
\be \Delta = a_1^6a_6 + a_1^5 a_3 a_4 + a_1^4 a_2 a_3^2 + a_1^4 a_4^2 + a_1^3 a_3^3 + a_3^4  \ee
On retrouve  bien l'expression du discriminant usuellement exprim\Ž \/ en fonction des coefficients $b_2, b_4, b_6$ , et le fait que le poids est $12$ (\cite {si} p 46).  \\ 
$\lozenge$
Notons maintenant $X\subset \mathbb A^{3m+2}$ l'ouvert des couples $(A,B)$ qui d\'efinissent  (par (2.3)) des courbes  hyperelliptiques  $C(A,B)$ non singuli\`eres, et soit $S  = \mathbb A^{3m+2} - X$ le ferm\'e    compl\'ementaire.  Nous allons d\'ecrire $S$, et donc $X$, en prouvant  que  $S$ est l'hypersurface  d'\'equation $\Delta \ne 0$, et qu'elle est irr\'eductible;  en fait on prouve mieux, la forme $\Delta$ est irr\'eductible.  Dans notre situation, caract\'eristique deux,  on notera que le discriminant usuel perd  une partie de sa signification, et d'ailleurs n'est plus irr\'eductible  comme on le voit imm\'ediatement avec le discriminant cubique\footnote {C'est un fait g\'en\'eral: en caract\'eristique deux le polyn\^ome de Vandermonde  est sym\'etrique, et   dans l'alg\`ebre des polyn\^omes sym\'etriques est  un polyn\^ome irr\' eductible. Le discriminant est le carr\Ž \/ du Vandermonde.}
$$ 27 {a_0}^2 {a_3}^2 + 4a_0 {a_2}^3 - 18a_0a_1a_2a_3 + 4{a_1}^3 a_3 - {a_1}^2 {a_2}^2$$
qui se r\'eduit \ˆ \/  $(a_0a_3+a_1a_2)^2$ en caract\'eristique deux. Il faut  dans le probl\`eme qui nous concerne lui substituer $\Delta$. Pour se convaincre de cela,  d\'ebutons par un lemme qui identifie $S$ avec $\Delta = 0$:
\begin {lemma}  Soit  toujours   $C(A,B)$ la courbe d\'efinie par la paire $(A,B)\in \Gamma ({\cal O} (m)) \oplus \Gamma ({\cal O} (2m))$.   Alors $C(A,B)$ est non singuli\`ere si et seulement si $\Delta (A,B) \ne 0$.  Si $\Delta (A,B) = 0$, et si le discriminant de la forme $A$ est non nul, la courbe $C(A,B)$ est nodale (eventuellement r\'eductible).\end{lemma}
\begin{proof} Noter que si, soit $A = 0$, soit $B = 0$, alors $C(A,B)$ est purement ins\'eparable sur $\mathbb P^1$ dans le premier cas, et r\'eductible dans le second. Par ailleurs 
$$ 0\times \mathbb A^{2m+1} \cup \mathbb A^{m+1} \times 0 \subset \{ \Delta = 0\} $$
On  peut donc se limiter  maintenant \ˆ \/ des couples $(A,B)$ avec $A\ne 0, B\ne 0$. 
Supposons en premier le couple $(A,B)$ tel que $a_0\ne 0$, ou bien $a_m\ne 0$.  Supposons par exemple    $a_0\ne 0$. 
Alors $A(1,0) \ne 0$,  en cons\'equence le rev\^etement $C(A,B) \to \mathbb P^1$ \'etant non ramifi\Ž \/ \ˆ \/ l'infini (Lemme 2.3), 
les points singuliers \'eventuels  de $C(A,B)$ sont au dessus de l'ouvert $Y\ne 0$. Si $P$ est un tel point,  au dessus du point de coordonn\'ees $(x,y)$, 
avec $y\ne 0$, alors 
$$A(x,y) = A'_X(x,y)^2B(x,y) + B'_X(x,y)^2 =  0$$
  ce qui implique $ a_0^2 \; \Delta (A,B) = 0$. Donc $\Delta (A,B) = 0$.  La r\'eciproque est claire:  si $\Delta (A,B) = 0$ alors  $A(X,Y) =   A'_X (X,Y) ^{2} B(X,Y) + B'_X (X,Y)^2 = 0 $ doit avoir une solution non triviale, et comme $a_0 \ne 0$, cette solution est $(x,y)$ avec $y\ne 0$. Ce point d\'efinit un point singulier de $C(A,B)$. 

L'argument est le m\^eme si $a_m \ne 0$. Cela prouve   en particulier  le r\'esultat  si $a_0\ne 0$, ou bien $a_m\ne 0$.  Notons en g\'en\'eral que si on effectue un changement lin\'eaire des variables, les courbes $C(A,B) $ et $C(A\circ \sigma , B\circ \sigma)$ sont isomorphes.  Si $a_0 = a_m = 0$, on se ram\`ene au cas $a_0 \ne 0$ en effectuant un tel   changement  de variables $x\to x, \, y\to \alpha x + y$, avec $A(1,\alpha) \ne 0$. Cela est toujours possible du fait que le corps de base est suppos\Ž \/ alg\'ebriquement clos. L'argument est licite du fait que  $\Delta (A,B)$ est semi-invariant (Lemme 3.4 (3)). Supposons les derni\`eres conditions   r\'ealis\'ees. On peut supposer que $a_0\ne 0$; soit un  point singulier $(z_0,x_0)$ au-dessus de $(x_0,1)$. Comme $A_x(x_0, 1) \ne 0$, il est clair que la courbe affine d'\'equation $z^2 + A(x,1)z  = B(x,1)$ qui a pour forme locale en $(z_0,x_0)$
$$(z -z_0)^2 + (z -z_0)(x - x_ˆ) A^*(x) = (x - x_0)^2 B^*(x),\,\, A^*(x_0) \ne 0$$
 pr\'esente un point double en $(z_0,x_0)$. Si $B = 0$, la courbe correspondante est r\'eductible.
  \end{proof}
Du fait que la caract\'eristique est deux,  on dispose aussi de l'action  du groupe additif $\mathbb G_a^{m+1}$ sur $\mathbb A^{3m+2}$ d\'ecrite ci-dessus.  Pour une autre description de l'hypersurface $\Delta = 0$, introduisons    le sous ensemble  $Z \subset \mathbb A^{3m+2}$  d\'efini par  la condition $ (A,B) \in Z  \Longleftrightarrow  \\ A\ne 0, B\ne 0,  $ et  il existe une forme lin\'eaire  $ \ell \ne 0 $   telle que 
\be  \ell / A\, ,  \,\ell^2  / B \ee
Il  a \'et\'e not\Ž \/  que $\Delta (\alpha . (A,B)) = \Delta (A,B)$,   par ailleurs il est clair que $X$ est stable par cette action, du fait que $C(\alpha .(A,B)) \cong C(A,B)$. 
\begin {thm} Le polyn\^ome $\Delta (A,B)$ est irr\'eductible, et    $X $ est l'ouvert  compl\'ementaire de l'hypersurface $\Delta \ne 0$.  On a 
\be \{ \Delta = 0\} =  \overline {\mathbb G_a^{m+1} . Z} \ee
\end{thm}
 \begin{proof}
 
  Il est imm\'ediat de voir  que si $(A,B) \in Z$ alors  la courbe  $C(A,B)$ est singuli\`ere.  Rappelons que  $S = \{ \Delta = 0\} $ est le compl\'ementaire de $X$,   lieu des paires $(A,B)$ telles que $C(A,B)$ est singuli\`ere (Lemme 3.5).  Observons d'abord que $Z$ est un ferm\Ž \/ irr\'eductible\footnote {Si $m=1$, les \'equations de $Z$ sont $b_1 = 0, \, a_0^2b_2 = a_1^2b_0.$}  de $(\mathbb A^{m+1} - 0) \times (\mathbb A^{2m+1} - 0)$ de  dimension $3m$. L'image de  $\pi (Z)$  de $Z$ dans $\mathbb P^m \times \mathbb P^{2m}$  qui est  aussi l'image du morphisme 
$$ \mathbb P^1\times \mathbb P^{m-1}\times \mathbb P^{2m-2} \longrightarrow, \,([l],[A^*],[B^*]) \mapsto ([lA^*],[l^2B^*])  $$
est un ferm\Ž \/ irr\'eductible. Comme $Z$ est stable par l'action de $\mathbb G_m^2$ (d\'efinie dans la section 3.2), le r\'esultat est clair, sauf peut \^etre la dimension. Mais  il est imm\'ediat de voir que $\dim \pi (Z) = 3m-2$, et donc $\dim Z = 3m$.
Montrons maintenant  que $S \subset \mathbb A^{3m+2}$  est d\'ecrit   par
\be  S = \overline {\mathbb G_a^{m+1} . Z}  \ee
On peut d'abord observer que l'adh\'erence de $Z$ dans $\mathbb A^{3m+2}$  contient les couples $(A,0)$. En effet, si on \'ecrit $A =  lA_0$ pour une certaine forme lin\'eaire $l$, alors si $B_0$ est une quelconque forme de degr\Ž \/ $2m-2$
$$\lim_{t\to 0} (A , t l^2 B_0) = (A,0)$$
Un argument analogue montre que $\overline Z \cap( 0\times \mathbb A^{2m+1}) $ est l'ensemble des couples $(0,B)$, $B$ \'etant de discriminant nul\footnote {On peut de mani\`ere analogue voir directement  que $0\times \mathbb A^{2m+1} \subset \overline {\mathbb G_a^{m+1} . Z}$. Il suffit de voir (si $m\geq 2$) que l'hypersurface d\Žfinie par l'annulation  du discriminant d'une forme de degr\Ž \/ $2m$ n'est pas stable par l'action du groupe additif $(\alpha,B) \mapsto B + \alpha^2$. Par exemple en notant que la forme $X^{2m} + a XY^{2m-1}$ n'a que des facteurs simples si $a\ne 0$.}. Observons maintenant que  l'adh\'erence  $\overline {\mathbb G_a^{m+1} . Z}$  est un ferm\Ž \/ irr\'eductible de codimension un. Soit  le morphisme $\varphi: \mathbb G_a^{m+1} \times Z \rightarrow \mathbb A^{3m+2}$:
$$\varphi (\alpha , (A,B)) = (A,F = B+\alpha A + \alpha^2)$$
Si $(A,F)$ est dans l'image de $\varphi$, la fibre en ce point est en bijection avec les $\alpha $ tels que $B = F + \alpha A + \alpha^2$ soit  de discriminant nul, donc est de dimension $m$. De la sorte  la dimension cherch\'ee est bien $3m+1$.\\
Pour conclure,  soit  $(A,B) \in S$, et montrons que $(A,B) \in  \overline {\mathbb G_a^{m+1} . Z} $. En vertu de ce qui vient d\^etre dit, il suffit de v\'erifier cela en dehors d'un ferm\Ž  \/ de codimension $\geq 2$, en particulier on peut supposer que $A\ne 0$, et $B\ne 0$. Supposons alors  que le point  $(x_0,y_0)$ soit un point "singulier" de $C(A,B)$,  avec par exemple $y_0 \ne 0$.   Alors  $A(x_0,y_0) = 0$, et quitte \ˆ \/  effectuer
une substitution $(A,B) \mapsto (A,B+\alpha A + \alpha^2)$, on peut s'arranger pour que  $B(x_0,y_0) = 0$.  Si $B = 0$, la conclusion est claire, car alors $(A,B) \in \mathbb G_a^{m+1} \overline Z$.  Dans le cas contraire  on doit avoir  $B'_X(x_0,y_0) = 0$, et alors ayant $A(x_0,y_0) = B(x_0,y_0) =  B'_X(x_0,y_0) = 0$, on obtient $(A,B) \in  Z$. L'inclusion oppos\'ee est claire.  En conclusion on a 
$$S = \{ \Delta = 0\} = \overline {\mathbb G_a^{m+1} . Z}$$
En particulier cela montre l'irr\'eductibilit\Ž \/ de l'hypersurface $S$.  

Pour  terminer la preuve du th\'eor\`eme 3.6 reste \ˆ \/ voir que  la forme $\Delta$ est irr\'eductible. Pour un r\'esultant  ou discriminant ordinaire en caract\'eristique z\'ero, c'est un fait bien connu \cite{gkz}. Du fait que $S$ est une hypersurface  irr\'eductible, il suffit de voir que $\Delta$ n'est pas  de la forme $\Delta_0 ^{e}$ pour un $e \geq 2$. 

Pour cela consid\'erons les formes g\'en\'eriques 
$$A = \prod_{i=1}^m (X- \alpha_i Y) \in k[\alpha_1,\cdots,\alpha_m] = k[\alpha] , \quad B = \sum_{j=0}^{2m} b_jX^{2m-j}Y^j \in k[b_0,\cdots,b_{2m}]$$
  $\alpha_1,\cdots,\alpha_m, b_0,\cdots,b_{2m}$ \'etant une collection d'ind\'etermin\'ees. On a sous ces conditions
$$\Delta (\alpha, b) =   \prod_{i=1}^m \,({A'_X}^2 B + {B'_X}^2) (\alpha_i,1))  $$
En d'autres termes, on se place dans 
$$R =  k[a_1,\cdots,a_m,b_0,\cdots,b_{2m}] \subset  S = k[\alpha_1,\cdots,\alpha_m,b_0,\cdots,b_{2m}]$$
les corps de fractions formant une extension galoisienne de groupe $\mathfrak S_m$.
Supposons avoir $\Delta = \Delta_0 ^{e}$, avec $\Delta_0$  n\'ecessairement irr\'eductible.   Posons
$$ G(X) = G(\alpha,b)(X) = ({A'_X}^2 B + { B'_X}^2) (X,1)  $$
 alors  le polyn\^ome $G$ a ses coefficients dans  $R$, et  par d\'efinition
$$\prod_{i=1}^m G(\alpha_i) = \Delta_0^{e}  $$
Observons que $G(\alpha_1) \in S$ est irr\'eductible. Il ne peut avoir de facteur irr\'eductible dans  $k[\alpha]$. Un tel facteur  en effet diviserait dans $k[\alpha]$ les polyn\^omes  d\'eduits de $G$ par les sp\'ecialisations $B = X^{2m}$, et $B = Y^{2m}$, imposant \ˆ \/  ce polyn\^ome d'\^etre un facteur de $A'_X (\alpha_1)$, situation clairement impossible. Du fait  que le degr\Ž \/ en les variables $b_0,\cdots,b_{2m}$ est  deux,  une d\'ecomposition en facteurs irr\'eductibles ne peut  donc \^etre que de la forme $G(\alpha_1) = uv \,(u,v\in S)$. On voit imm\'ediatement  comme cons\'equence de l'argument pr\'ec\'edent que $u$ et $v$ doivent alors \^etre de degr\Ž \/ un en les variables $b_i$.  Posons
$$u = \phi_0 + \sum_{i=0}^{2m} u_ib_i, \,\, v = \psi_0 + \sum_{j=0}^{2m} v_jb_j\quad (\phi_0,\psi_0,u_i,v_j \in k[\alpha_1,\cdots,\alpha_m]  $$ 
En developpant  le produit $uv$ (2.17),  et en identifiant avec $G(\alpha_1)$, on trouve  $\phi_0\psi_0 = 0$, pour tout $k=0,\cdots,2m$
$$u_kv_k =  \begin{cases} \alpha_1^{2(2m-k-1)} &\text {si $k$ est impair }\\  0 & \text {si $k$ est pair}.\end{cases}  $$
et $u_iv_j + u_jv_i = 0 $ si $i\ne j$. Supposons par exemple $\phi_0 = 0$; alors en utilisant les relations qui pr\'ec\`edent on trouve que pour tout $k$, $u_k \ne 0$, ceci du fait que $\psi_0 u_k = {A'_X}^2(\alpha_1){\alpha_1}^{2m-k}$.  Cela implique   ais\'ement $v_k = 0$ pour tout $k$, ce qui est impossible. Ainsi $G(\alpha_1)$ est irr\'eductible dans $k[\alpha,b]$. Il est imm\'ediat d'en d\'eduire que la norme $\prod_i G(\alpha_i)$ est un \'el\'ement irr\'eductible de l'alg\`ebre de polyn\^omes $k[a,b]$.

Cela montre finalement que $e = 1$, donc que $\Delta (A,B)$ est irr\'eductible.  Ce raisonnement permet de prouver de nouveau, sans r\'ef\'erence \ˆ \/ l'argument  qui pr\'ec\`ede l'irr\'eductibilit\Ž \/ de $\Delta$.
  \end{proof}
 
\sectiono { La stratification de Hasse du champ  ${\cal H}_g $ }
\subsection{${\cal H}_g$ comme champ quotient}

 Dans  cette section nous prouvons notre r\'esultat principal qui est la d\'escriptionde  ${\cal H}_g$ comme un champ quotient. Comme dans les sections ant\'erieures, le champ ${\cal H}_g$ est suppos\'e d\'efini sur le corps alg\'ebriquement clos $k$ de carac t\Ž ristique deux. Comme   indiqu\Ž \/ dans l'introduction la description que nous allons proposer peut  s'interpr\'eter  comme une variante de l'utilisation de la forme de Weierstrass (ou de Rosenhaim) g\'en\'eralis\'ee \cite {i}.  On observera cependant que si $g=1$, on a ${\cal H}_1 \ne {\cal M}_{1,1}$. 

Soit une courbe hyperelliptique $p: C\to S$, et consid\'erons le rev\^etement s\'eparable de degr\Ž \/ deux $\pi:C\to D = C/\tau$. On garde les notations de la section 2, en particulier $\tau$ est l'involution hyperelliptique.  On d\'eduit de cette donn\'ee la  ${\cal O}_D$  alg\`ebre localement libre ${\cal E} = \pi_*({\cal O}_C)$, et l'extension  correspondante (3.2). Localement  pour la topologie \'etale  sur $S$, cette donn\'ee se trivialise en  (Lemme 3.2)
\be  D = \mathbb P^1_S,\quad {\cal E} = {\cal O}_D \oplus {\cal O}_D(-m)  \ee
La structure d'alg\`ebre sur ${\cal E}$ est  alors compl\`etement d\'ecrite par la donn\'ee d'un couple de sections  (\S 3.4)
$$ (A,B) \in \Gamma ({\cal O}_D(m)) \times \Gamma ({\cal O}_D(2m))  $$
conjointement avec le sous-fibr\Ž \/ ${\cal O}_D$, engendr\Ž \/ par l'\'el\'ement unit\Ž \/.   Il a \'et\'e\/ rappel\Ž \/ (Proposition 2.2) que l'ensemble des trivialisations  de  la donn\'ee $(D , {\cal E}, \tau)$, c'est  \ˆ \/  dire le $S$-sch\'ema
\be \Isom \left( (\mathbb P^1_S,{\cal O} \oplus {\cal O}(-m)) , (D,{\cal E})\right )  \ee
est un torseur sous le groupe des automorphismes du mod\`ele,   et r\'eciproquement  tout torseur sous ce groupe est de cette forme.  Soit dans ce contexte le groupe $G\times S$,  donn\Ž \/ comme  $G = \Aut \left( (\mathbb P^1_S,{\cal O} \oplus {\cal O}(-m) \right)$\footnote {Rappelons que les automorphismes sont assujettis \ˆ \/ fixer la section unit\Ž \/.}.  On a vu  que ce groupe  se r\'ealise comme un sous-groupe  du groupe qui lin\'earise universellement  le fibr\Ž \/ mod\`ele ${\cal O}_{\mathbb P^1}  \oplus {\cal O}_{\mathbb P^1}(-m)$.  Il s'ins\`ere   dans ce cas particulier dans une extension
\be  1 \rightarrow H \rightarrow G \rightarrow  \PGL (1)  \rightarrow 1  \ee
 $H = \Aut \left(  {\cal O}_{\mathbb P^1}  \oplus {\cal O}_{\mathbb P^1}(-m) \right)$ \'etant  le groupe des automorphismes \'egaux \ˆ \/ l'identit\Ž \/ sur ${\cal O}_{\mathbb P^1}$\footnote { Dans cette situation il faut pr\'eciser que l'objet marque est un triplet ($\mathbb P^1, V, {\cal O}_{\mathbb P^1} \hookrightarrow V$). Les automorphismes sont ceux qui se r\'eduisent \ˆ \/ l'identit\Ž \/ sur la section unit\Ž \/}, c'est \ˆ \/ dire
$$H \; = \; \left ( \begin{pmatrix} 1 & \alpha\\ 0&\beta\\  \end{pmatrix} , \alpha \in \Gamma ({\cal O}_{\mathbb P^1} (m)), \,\beta \in \Gamma ({\cal O}_S^*) \right) $$
Sa structure pr\'ecise, cas particulier du lemme 2.3, est rappel\'ee  pour m\'emoire:
\begin {lemma} \begin{enumerate} \item Si $m$ est pair, $G \cong H\rtimes \PGL (1)$ (la suite ({\rm 4.3}) est scind\'ee). 
\item  Si $m$ est impair $G \cong \Gamma ({\cal O} (m)) \rtimes \GL (2) / \mu_m$.  \\
Dans l'assertion $(2)$ le produit semi-direct est d\'efini relativement \ˆ \/ l'action \'evidente de  $\GL (2)/\mu_m$ sur $\Gamma ({\cal O}(m))$. On a  donc dans tous les cas $G \cong \Gamma ({\cal O} (m)) \rtimes \GL (2) / \mu_m$. \end{enumerate}\end{lemma}
\begin{proof}1) On consid\`ere pour tout entier $n$,  le sous-groupe  du groupe des matrices diagonales $G_m \subset  \GL (2)$, noyau de   $G_m \to G_m, \,\, \lambda \to \lambda^n$.
Ce sous-groupe isomorphe \ˆ \/ $\mu_n$, est not\Ž \/ $\mu_n$. Noter que si $n$ est pair, il n'est pas \'etale; cependant le groupe dual est toujours $\frac {\mathbb Z } {n\mathbb Z}$.  Le faisceau ${\cal O}(-m)$ admet une $\PGL (1)$-lin\'earisation canonique, qui provient par passage au quotient  et puissance $-m $ de  la lin\'earisation de ${\cal O}(1)$ de groupe $\GL (2)$.  Si cette derni\`ere est $\alpha_g: {\cal O} (1) \cong [g]^* ({\cal O} (1))$, la section de  (4.3) est 
$$[g] \mapsto ([g] , \alpha_g^{\otimes -m})$$
Dans le cas (2),  la lin\'earisation de ${\cal O} (m)$ sous $\GL (2)$ factorise par $\GL (2) / \mu_m$. Il en r\'esulte une section $\GL (2) / \mu_m \to G$; explicitement $\overline g \mapsto ([g],\alpha_g^{\otimes -m}$. Il est visible que l'image  du sous-groupe des homoth\'eties de $\GL (2)$ est le sous-groupe 
$$\begin{pmatrix} 1&0\\0&\star \\ \end{pmatrix} \subset H$$
On  en tire   le fait que $G =  \Gamma ({\cal O} (m)) \rtimes \GL (2) / \mu_m$.  Comme cet argument ne d\'epend pas de la parit\Ž \/ de $m$, le r\'esultat est donc valable dans les deux cas.  \end{proof}
 
Notons que le groupe $G$ agit naturellement sur l'espace affine $\mathbb A^{3m+2}$,  de sorte que l'ouvert $X$ est  $G$-stable.  Cette action m\'elange l'action de $\GL (2)$, plus pr\'ecis\'ement de   $\GL (2) / \mu_m$ d'une part  et celle de $H$. L'action de $H$ (exprim\'ee \ˆ \/ droite)  qui traduit les changements de trivialisations du fibr\Ž \/ ${\cal E}$ est
\be \begin{pmatrix} 1&\alpha \\ 0 & \beta \end{pmatrix} \begin{pmatrix} B\\A \\ \end{pmatrix} \;=\; \begin{pmatrix} \beta^2 B + \beta \alpha A + \alpha^2 \\ \beta A \\ \end{pmatrix}   \ee
Le couple de sections $(A,B)$ qui d\'ecrit la structure d'alg\`ebre sur ${\cal E}$  une fois trivialis\Ž \/,  s'interpr\`ete alors comme un morphisme $G$-\'equivariant  
$$\Isom \left( (\mathbb P^1_S,{\cal O} \oplus {\cal O}(-m)) , (D,{\cal E})\right )  \longrightarrow X $$
D\`es lors, on peut \'enoncer  le r\'esultat, en acceptant l'existence de  ${\cal H}_{g,\mathbb Z}$ \cite {rom} (pour la d\'efinition d'un champ quotient, voir par exemple \cite {lmb}):
\begin {thm} On a un isomorphisme de champs ${\cal H}_g  = {\cal H}_{g,k} \cong  [X/G]$, $[X/G]$ \'etant le champ quotient de $X$ par $G$.  En particulier ${\cal H}_g $ est   lisse.  En conclusion, et de  mani\`ere plus g\'en\'erale,  ${\cal H}_{g,\mathbb Z}$  est  lisse sur $\mathbb Z$.\end{thm}
 \qed
 \subsection{La stratification de Hasse}

Soit une courbe hyperelliptique d\'efinie sur le corps $k$.  On sait (Lemme 3.3) que le diviseur des points fixes  $F$ de l'involution $\tau$  a   dans la description  qui pr\'ec\`ede pour \'equation $A = 0$. Notons que  le sens de ce diviseur, qui est le diviseur de Weierstrass \cite{kl},  devient plus  clair si on fait intervenir la ramification sup\'erieure, c'est \`a dire  le conducteur de Hasse en chaque point fixe $P_i, \; (1\leq i \leq b) $  \cite {bm},\cite {jps}. Rappelons que si $P$ est un point fixe de $\tau$, et si $t\in \cal M_P$ est une uniformisante en $P$, le conducteur de Hasse $m_P$ est d\'efini par $v_P(\tau (t) - t) = m_P+1$.  On sait que $m_i$ est impair. En r\'esum\Ž \/
\begin{lemma} On a 
\be  F = \sum_{i=1}^b (m_i+1 ) P_i , \quad \deg (F) =  2m = 2g+2 \ee
  $m_i$  d\'esignant  le conducteur de Hasse  en $P_i$.  De plus si $C$ est une courbe  de  base $S$, alors $F$ est un  diviseur relatif sur $S$ de degr\Ž \/ $ 2m$.  Il existe un diviseur  $Z \subset D$, de degr\Ž \/  $m$, tel que $F = \pi^* (Z)$.\end{lemma} 
 D'une autre mani\`ere consid\'erons  les groupes de ramification sup\'erieurs; en $P_i$, le conducteur \'etant $m_i$,  il sont d\'efinis par  
$$\Gamma = (1,\tau) = \Gamma_0 = \cdots = \Gamma_{m_i} \varsupsetneqq \Gamma_{m_i+1} = 1$$
 Si on \'ecrit les multiplicit\'es $m_i+1\over 2$ de mani\`ere d\'ecroissante, on  d\'efinit  de la sorte 
une partition de $m$. Noter qu'\ˆ \/  toute partition $\mu = (\mu_1 \geq \cdots \geq \mu_b > 0)$ de $m$ on peut associer de cette mani\`ere  une  courbe hyperelliptique de genre $g$. On notera   pour une partition donn\'ee $\mu$,  ${\cal H}_g (\mu)$ le  champ des courbes hyperelliptiques a ramification fix\'ee de type $\mu$. On a le r\Ž sultat ais\Ž \/ suivant:
\begin{propo}
Pour toute  partition $\mu$ de $m = g+1$,  de longueur $b$, le champ ${\cal H}_g(\mu)$ est un sous champ localement ferm\Ž  \/ et lisse de   codimension  $g-b+1$  de ${\cal H}_g$. \end{propo} 

\begin{proof}  Il est clair en vertu du th\'eor\`eme 4.2, et du Lemme 4.3,  que ${\cal H}_g(\mu)$ est le champ quotient  $[X(\mu)/G]$, avec $X(\mu)$ le sous-sch\'ema  localement ferm\Ž \/ de $X$  dont les points sont les couples $(A,B)$, la forme $A$  \Žtant telle que $\Div (A)$ est de type $\mu$ (Lemme 4.3).  Le r\'esultat  suit.  \end{proof} 
On notera que si $\mu = (2,1,\cdots,1), \; b= m$, alors l'adh\'erence de   cette strate est un diviseur de ${\cal H}_g$. Pour les autres strates la codimension est  au moins deux. Par ailleurs,  la  formule de Crew pour le  calcul du $p$-rang  $r_C$  d'un rev\^etement  galoisien $\pi: C\to D$, de groupe $Q$,  g\'en\'eriquement \'etale de courbes propres, lisses et connexes sur un corps alg\'ebriquement clos de caract\'eristique $p > 0$ (voir par exemple \cite {ray}) 
\be   r_C - 1 =  \vert Q\vert (r_D - 1) + \sum_{y\in D} (e_y - 1) \ee
les $e_y$ d\'esignant les indices de ramification, montre que  le long d'une strate  ${\cal H}_g (\mu)$,  le $2$-rang est constant \'egal \ˆ \/ $b-1$. En particulier la strate ouverte est form\'ee des courbes ordinaires, et la strate ferm\'ee, de dimension $m = g-1,$ est constitu\'ee de courbes de $2$-rang z\'ero.
\begin {rem}  \end{rem} Pr\'ecisons le lien si $g=1$ entre les deux champs ${\cal H}_1$ et ${\cal M}_{1,1}$ (voir aussi \cite {av}). Notons d'abord  qu'il y a un  morphisme  naturel
\be  \psi: {\cal M}_{1,1} \longrightarrow {\cal H}_1 \ee
 donn\Ž \/ sur les objets par $\psi (C\to S, O) =  (C\to S, \tau)$, avec pour $\tau$ l'involution $x\mapsto -x$; l'oppos\Ž \/  $- x$ est  d\'efini relativement \ˆ \/   la loi de groupe sur $E$, d'unit\Ž \/ la section $O$. Comme cette involution est la seule fixant le  $S$-point $O$,  cela d\'efinit bien un foncteur, i.e. un morphisme de champs.  Montrons que ce morphisme est repr\'esentable,  plat,  fini de degr\Ž \/ quatre.  De mani\`ere plus pr\'ecise c'est le diviseur de Cartier universel ${\cal F}  \to {\cal H}_g$ (Lemme 3.3). Cela est \ˆ \/ pr\`es imm\'ediat. Soit  $S \to {\cal H}_1$, un point  d\'efini par la courbe $(C/S,\tau)$. Une section au dessus de $T$, du  $2$-produit fibr\Ž \/ ${\cal M}_{1,1}\times_{{\cal H}_1} S$ est la donn\'ee  d'un morphisme $T\to S$, d'une courbe elliptique $E \to T$, avec section nulle $O: T\to E$, et d'un isomorphisme \'equivariant relativement aux involutions $\phi: E\stackrel\sim \rightarrow C\times_S T$. La section $O$ d\'efinit   ainsi un $T$-point de $F\to S$. Il est clair que ${\cal M}_{1,1}\times_{{\cal H}_1} S \cong \cal F$.
 $\lozenge$
\sectiono {Le groupe de Picard de ${\cal H}_g$ }
 
Rappelons que si on regarde ${\cal H}_g$ comme \Žtant  d\Žfini sur un corps $k$ de caract\Ž ristique\footnote {plus pr\'ecis\'ement $p$ premier \ˆ \/ 2 et $g+1$.}   $p\ne 2$, alors Arsie et Vistoli (\cite {av}, thm 5.1, remark 5.5) ont  montr\Ž \/ que le groupe de Picard de ${\cal H}_g$ (voir \cite {av}, \cite {mum}, \cite {vg}  pour la d\'efinition du groupe de Picard) est 
\be \Pic ({\cal H}_g) = \begin{cases} \frac {\mathbb Z} { (8g+4) \mathbb Z} \;\; \text {si }\;$g$\; \text { est \;impair } \\
\frac {\mathbb Z} { (4g+2) \mathbb Z} \;\; \text {si} \; $g$ \; \text {est \;  pair} \end{cases}  \ee
On va prouver que ce r\'esultat subsiste sans restriction sur la caract\'eristique, donc si la caract\'eristique est deux, bien que la description du champ soit assez diff\'erente.  La d\'etermination du groupe $\Pic ({\cal H})$  lorsque ${\cal H} = [X/G]$  est un champ quotient d'une vari\'et\Ž \/ lisse $X$, est dans certains cas ais\'ee.  De la d\'efinition du groupe de Picard, on tire en effet
\be \Pic ({\cal H}) = \Pic_G (X)  \ee
o\ \/  $\Pic_G (-)$ d\Žsigne le groupe des faisceaux inversibles $G$-lin\'earis\'es. Le lemme suivant permet dans certains cas, et pour des $G$-vari\'et\'es affines,  d'expliciter le groupe de Picard \'equivariant.  Rappelons  que si $X$ est une $k$-vari\'et\Ž \/ normale (de type fini), on pose (d\'efinition de Rosenlich)
$U_k(X) = \Gamma (X , {\cal O}_X^*)/k^*$. On sait que ce groupe ab\'elien est  libre de type fini, et que  
$$U_k (X \times Y) = U_k (X) \times U_k (Y)$$
Si $G$ est un $k$-groupe alg\'ebrique affine lisse, alors $U_k (G) = \hat G$, le groupe des caract\`eres.
  \begin{lemma}  Soit $X = \Spec (R)$ le spectre d'une   $k$-alg\`ebre de type fini normale,  d\'efinie sur un corps alg\'ebriquement clos $k$. On suppose que $X$  est muni d'une action  r\'eguli\`ere d'un groupe alg\'ebrique lisse connexe $G$.   On a  alors une suite exacte\footnote {Le morphisme $\beta$ traduit l'obstruction \ˆ \/ ce qu'un faisceau inversible $G$-invariant soit $G$-lin\'earisable; il ne sera pas utilis\Ž \/.}:
  \be 0 \rightarrow  \frac {\hat G } { \langle \chi \rangle} \rightarrow \Pic _G (X)\stackrel\alpha  \rightarrow \Pic (X)^G  \stackrel\beta  \rightarrow  H^2(G,k^*)  \ee
dans laquelle le terme $\Pic (X)^G$ d\'esigne le sous-groupe des  \'el\'ements $G$-invariants de $\Pic (X)$, et $\langle \chi \rangle$ d\'esigne le sous-groupe engendr\Ž \/ par les ''poids'' des \'el\'ements de $U_k (X)$. En particulier si $\Pic (X) = 0$ ($R$ est factoriel), on a $\Pic_G (X) = \frac {\hat G} { \langle \chi \rangle}$.
 \end{lemma} 
  \begin{proof} Pr\'ecisons que dans la suite  (5.3), le morphisme $\alpha$ est l'oubli de l'action de $G$.    Le groupe $G$ \'etant un groupe alg\'ebrique  connexe (lisse ou non),   agissant sur une vari\'et\Ž \/ normale  $X$,   on sait qu'il y a une suite exacte
 \be 0\rightarrow H^1(G,\Gamma (X,\Cal O_X)^*) \rightarrow \Pic (X)^G\stackrel\alpha  \rightarrow \Pic (X) \rightarrow \Pic (G) \ee
La description du noyau de $\alpha$ traduit le fait qu'une lin\'earisation du fibr\Ž \/ trivial de rang un $X\times \mathbb A^1 \to \mathbb A^1$ est  d\'ecrite  au niveau des points par
$$(g , x , t) \mapsto (gx , \varphi (g,x)t)$$
la fonction $\varphi$ v\'erifiant la relation de cocycle 
$$\varphi (gg' , x) \,= \,\varphi (g,g'x) \,\varphi (g',x)$$
 et \'etant d\'etermin\'ee modulo un cobord $(g,x)\mapsto {f(gx)\over f(x)}$. 
Dans la situation qui nous occupe $X$ est un ouvert affine d'un espace affine, on a donc $\Pic (X) = 0$,   on est ainsi essentiellement ramen\Ž \/ au calcul du groupe     $ H^1(G,\Gamma (X,{\cal O}_X)^*) $. Avec les hypoth\`eses faites sur $X$,  le r\'esultat de Rosenlich rappel\Ž \/ au-dessus montre qu'un cocycle est repr\'esent\Ž \/ par un caract\`ere de $G$. Le reste en d\'ecoule  imm\'ediatement.  \end{proof} 
   
 On est  maintenant en mesure de calculer le groupe de Picard du champ ${\cal H}_g$ sur un corps alg\'ebriquement clos $k$ de caract\'eristique deux, et prouver le r\'esultat annonc\'e dans l'introduction. La conclusion qui peut surprendre,  est que le r\'esultat est le m\^{e}me qu'en caract\'eristique $p\ne 2$. Cela est sans doute \`a rapprocher du fait que pour les rev\^etements doubles la lissit\'e est pr\'eserv\'ee.  On notera cependant  que dans cette description la premi\`ere classe de Chern du fibr\Ž \/ de Hodge n'est pas  dans tous les cas un g\'en\'erateur  du groupe de Picard (Remarque 4.5).  Tenant compte de {\rm  \cite {av}}, le r\'esultat suivant montre  (5.1)  est vrai  en toute caract\'eristique.

 \begin{thm}  Le groupe de Picard du champ ${\cal H}_g = {\cal H}_{g,k} \;(g\geq 1)$, suppos\Ž \/ d\'efini au-dessus  d'un corps alg\Žbriquement clos $k$ de caract\Ž ristique $p=2$,  est : 
 \be \Pic ({\cal H}_g) = \begin{cases} \frac {\mathbb Z} { (8g+4) \mathbb Z} \;\; \text {si }\;$g$\; \text { est \;impair } \\
\frac {\mathbb Z} { (4g+2) \mathbb Z}  \text\;\; {si} \; $g$ \; \text {est \;  pair} \end{cases}  \ee
 \end{thm} 
 \begin{proof} 
 Le  th\'eor\`eme 4.2  dit que ${\cal H}_g = [X/G]$, le groupe $G$ \'etant  d\'ecrit  dans le lemme 4.1, et  avec $X = \mathbb A^{3m+2} - \{ \Delta = 0\}$. Il a \'et\Ž \/ observ\Ž  \/  d'autre part  ($\S 3$) que le diviseur irr\'eductible et r\'eduit  $\Delta = 0$ est invariant sous l'action naturelle de $G$, et que plus pr\'ecis\'ement  
 $\Delta $ est un polyn\^ome irr\'eductible semi-invariant. Si $Ê\chi \in \hat G$ est le  poids de $\Delta$,  alors (Lemme 5.1) on peut conclure:
 $$\Pic ({\cal H}_g) = \frac  {\hat G} {\mathbb Z \chi}  $$
   Comme le groupe $G$ est  produit semi-direct  de $\GL (2) / \mu_m$ par un groupe unipotent,  il en  d\'ecoule que $\hat G =  \widehat {\GL (2)/\mu_m}$ et donc  $\hat G =   \mathbb Z \psi$ avec  $\psi = \det^{m}$ si $m$ est impair, sinon $\psi = \det^{\frac m 2} $ si $m$ est pair. Il reste donc \ˆ \/  trouver le poids $\chi = \psi^{e}$ de $\Delta$.  Il est suffisant  pour cela de faire agir le sous-groupe des matrices diagonales, c'est \ˆ \/ dire $G_m$ via  $\lambda \mapsto \begin{pmatrix} \lambda & 0 \\ 0& \lambda \end{pmatrix} $.  Alors $A$ est de poids $m$ (relativement \ˆ \/  l'action de $G_m$)  , et ${A'_X}^2 B + {B'_X}^2$ est de poids $4m-2$. Donc le poids  correspondant de $\Res (A , {A'_X}^2 B + {B'_X}^2)$ est $2m(4m - 2)$.  Finalement celui de $\Delta (A,B)$ est $2m(4m  -2) - 2m = 2m(4m - 2)$ . Si on revient  au calcul de l'exposant $e$, on trouve en conclusion 
 $$e = \begin{cases} 4m-2 \;\; \text { si } m \;\text {impair }\\
 8m-4 \;\;\text {si}\;\; m \;\text {pair} \end{cases}$$ 
  ce qui est  le r\'esultat annonc\Ž.  \end{proof}     
 \begin{rem}    \end{rem}
  Si $g=1$,  le fait que ${\cal M}_{1,1} \ne {\cal H}_1$ (Remarque 4.5),  et le fait que ces deux champs ont  un groupe de Picard identique,  admet une  explication simple.  On peut aussi  d\'eduire du th\'eor\`eme 5.2 le  r\'esultat  fameux  $\Pic ({\cal M}_{1,1}) = \frac {\mathbb Z} {12 \mathbb Z}$. Pour que les choses soient claires, on notera que le champ ${\cal M}_{1,1}$, \ˆ \/ la diff\'erence de ${\cal H}_1$, n'est pas un champ quotient \cite {bck}.  Ceci \'etant  il n'est pas difficile de voir que $\psi^\star: \Pic ({\cal H}_1) \to \Pic ({\cal M}_{1,1})$ est un isomorphisme. Cela r\'esulte par exemple de l'interpretation de $\psi$, donn\'ee dans la remarque 4.5. En effet,  $X$ ayant  la signification  des sections 3 et 4,  et donc fournissant un atlas  $X\to {\cal H}_1$, on a remarqu\Ž \/ (loc.cit.)  que le carr\Ž  \/  suivant est $2$-cart\'esien:
  $$\CD {{\cal M}_{1,1}} @>\psi>> {{\cal H}_{1,1}}  \\
  @VVV @VVV\\
  \cal F @>\pi >> X \endCD  $$
  Ainsi  un faisceau inversible  sur ${\cal M}_{1,1}$ est incarn\Ž \/ par un faisceau inversible sur $\cal F$. Si ${\cal F}' = \cal F\times_X \cal F$, il est facile de voir que ce faisceau inversible  se trouve \^etre muni d'une donn\'ee de descente  relativement \ˆ \/ $\psi$.  On utilise pour cela la translation qui relie les deux sections  au dessus de ${\cal F}'$.  Ainsi il provient d'un unique faisceau inversible sur ${\cal H}_1$. 
$\lozenge$

Revenons au cas g\'en\'eral.  Nous allons identifier dans $\Pic ({\cal H}_g) = \frac {\mathbb Z} {e\mathbb Z}$ ($e$ \'etant comme dans le th\'eor\`eme 5.2)  les   deux classes naturelles
\be \lambda = c_1 (\mathbb E) = \det (\R p_! (\omega_{C/S}))\ee
 le d\'eterminant du fibr\Ž \/ de Hodge  \cite {av}, \cite {br} et la classe du diviseur ${\cal H}_g(2,1,\cdots,1)$ param\'etrant les courbes non ordinaires. On notera $\left [{\cal H}_g(2,1,\cdots,1)\right ]$ cette classe.
 \begin{propo} Dans $\Pic ({\cal H}_g) = \frac {\mathbb Z} {e\mathbb Z}$ avec $e = 4m-2$ si $m$ impair, sinon $e = 8m-4$, on a:\begin{enumerate}
 \item   $\lambda =  \begin{cases} {m-1\over 2}\;\; \text { si } m \;\text {impair }\\
 m-1 \;\; \text { si } m \;\text {pair }\\\end{cases} $ 
 \item     $\left [{\cal H}_g(2,1,\cdots,1)\right ] =  2\lambda$.\end{enumerate} 
 En particulier $\lambda$ est g\'en\'erateur si $m$ est pair, ou bien $m \equiv 3 \pmod 4$. \end{propo} 
 \begin{proof}  (1) Soit $p: C\to S$ une courbe hyperelliptique de base $S$, $q: D\to S$ le fibr\Ž \/ en coniques quotient de $C$ par l'involution $\tau$, et $\pi: C\to D$ le rev\^etement de degr\Ž \/ deux correspondant. Le faisceau inversible $\lambda \in \Pic ({\cal H}_g)$ est d\'efini par 
 $$\lambda (C\to S) =  \det (\R p_!  (\omega_{C/S})) = \det  (p_\star (\omega_{C/S} )) $$
 vu que $R^1p_* (\omega_{C/S} = {\cal O}_S$. Si $R\subset C$ est le diviseur de ramification relatif,  la formule de ramification donne 
 \be \omega_{C/S} = \pi^* (\omega_{D/S}) \otimes  {\cal  O} (R) \ee
 Par ailleurs on sait que $R = F$,  le diviseur des points fixes de l'involution hyperelliptique $\tau$, et  que ${\cal O} (R) = \pi^* ({\cal L}^{-1})$ (Lemme 3.3).  Donc finalement 
 $$ \lambda (C\to S) = \det  ( q_\star  \pi_\star (\pi^* (\omega_{D/S} \otimes {\cal L}^{-1}))$$
 ce qui,  par la formule d'adjonction, et la suite exacte (3.2), conduit \ˆ \/
$$ \lambda (C\to S)  = \det \R p_\star (\omega_{C/S}) = \det q_\star (\omega_{D/S}\otimes {\cal L}^{-1})  $$
 
Pour identifier ce faisceau inversible, on se place sur l'atlas $X$ d\'ecrit dans la section 3, alors $D = \mathbb P^1 \times X$  et $ {\cal L} = {\cal O}(-m)$. L'expression (5.6) d\'efinit  un module libre de rang un sur $k[X]$, \'equip\Ž \/ d'une action de $G$,  donc d\'efinie par un caract\`ere de $G$. Ce caract\`ere est facile \ˆ \/ expliciter\footnote { On peut de mani\`ere plus g\'en\'erale d\'ecrire le fibr\Ž \/ de Hodge $\mathbb E$. Il correspond  au $(k[X],G)$-module libre de rang $g$,  $\Gamma ({\mathbb P^1}_X  , \omega_{\mathbb P^1_X} \otimes {\cal O} (m))$.  Le r\'esultat est le module libre ${\cal O} (X)[X,Y]_{m-2}$ des formes  de degr\Ž \/ $m-2$ \ˆ \/ coefficients dans l'anneau des fonctions ${\cal O}(X)$. Ce module est  muni de l'action de $G$, $g.F = \det (g) F.g^{-1}$ qui d\'ecoule de la lin\'earisation produit tensoriel.}. Il faut  pour cela pr\'eciser la $G$-lin\'earisation support\'ee par  $\omega_{D/S}\otimes {\cal L}^{-1}$. C'est le produit tensoriel d'une part de la $G$-lin\'earisation  de $\omega_{D/X)} = {\cal O}_D(-2)$  qui provient de celle  canonique,  sous $\PGL (1)$, donc avec les notations de la section 3,  d\'ecrite par $g\in \GL (2) \mapsto  \det (g) \alpha_g^{\otimes (-2)}$. Pour  le second facteur ${\cal L}^{-1}$, la $G$-lin\Ž arisation est   donn\'ee par $g\in \GL (2) \mapsto  \alpha_g^{\otimes m}$.  En conclusion, on trouve que la $G$-lin\'earisation  de $\omega_{D/X)} \otimes {\cal L^{-1}}$ est  d\'efinie  par 
$$g\in \GL (2) \mapsto \det g \;\alpha_g^{\otimes (m-2)}  $$
En particulier la matrice diagonale $\mu 1_2$ agit par le facteur $\mu^{m(m-1)}$; en comparant au g\'en\'erateur $\psi$ (Th\'eor\`eme 5.2), le r\'esultat  (1) en d\'ecoule.\\
 (2)  Le discriminant $\Delta_m (A)$ d'une forme binaire $A$ de degr\Ž \/ $m$, \'etant  d\Ž fini  par
$$\Res_{m,m-1} (A,A'_X) = a_0 \Delta_m (A) $$
L'\'equation  du diviseur ${\cal H}(2,1\cdots,1)$ est, au niveau de l'atlas $X$ donn\'ee par  $\Delta_m (A) = 0$.    Le discriminant est  une forme de degr\Ž \/ $2m-2$ en les $a_i$, de sorte que  relativement au sous-groupe  \ˆ \/ un param\`etre $\mu \mapsto \mu.1_2$, son poids est $m(2m-2)$. Le r\'esultat en d\'ecoule  du fait du th\'eor\`eme 4.2. 
\end{proof} 
\begin{rem}\end{rem}
Le r\'esultat explicite de dessus montre que $\lambda = c_1(\mathbb E)$ engendre le groupe cyclique $\Pic ({\cal H}_g)$ si et seulement si $m $ est pair, ou bien $m \equiv 3 \pmod {4}$; sinon $\lambda$ est d'ordre deux.  C'est le cas si $g = 2$. Si $g=2$, et en caract\'eristique $\ne 2,3$, Vistoli  \cite {vg} a montr\Ž \/ que le l'anneau de Chow de ${\cal H}_2$ sur $\mathbb Z$  est d\'ecrit par 
$$A  ({\cal H}_2) = \mathbb Z[\lambda_1 , \lambda_2] / (10\lambda_1 , 2\lambda_1^2 - 24 \lambda_2)  $$
o\ \/ $\lambda_i = c_i (\mathbb E), \,\, i=1,2$. Il serait int\'eressant d'\'etendre  cette description au cas $p = 2$.  

On peut se demander si la propri\'et\'e de lissit\'e subsiste en les points du bord de la compactification stable   $\overline {{\cal H}}_g$.  Rappelons qu'une courbe hyperelliptique de genre $g \geq 2$ stable est donn\'ee par une courbe stable $C$ de genre $g$, munie d'une involution $\tau$ telle que $C/\tau$ soit de genre z\'ero, cette involution est alors unique.  Dans ce contexte, le probl\`eme propre \`a la caract\'eristique deux, vient du fait qu'on ne peut \'eliminer l'eventualit\'e de points g\'en\'eriques de composantes, des $\mathbb P^1$,  fix\'es par $\tau$. Cela appara\^{i}t en codimension $\geq 2$.  En ces points l'anneau de d\'eformation universelle n'est pas lisse.  Il serait utile dans ce cas de  pr\'eciser   la nature des obstructions \`a la d\'eformation des courbes hyperelliptiques stables, c'est \`a dire  d\'ecrire les \'equations  qui d\'efinissent l'anneau de la d\'eformation universelle. Ce probl\`eme avec des r\'eponses partielles est abord\'e dans \cite{bmg}.
$\lozenge$

  \end{document}